\documentclass[12pt]{amsart} 

\usepackage[colorlinks=true, linkcolor=black ,  citecolor=blue, urlcolor=black]{hyperref}
\usepackage{cleveref} 
\usepackage[utf8]{inputenc}

\title{Homotopy Techniques for Analytic Combinatorics in Several Variables}

\author{Kisun Lee }
\address{Department of Mathematics, University of California San Diego, 9500 Gilman Drive, La Jolla, CA, USA 92093}
\email {kisunlee@ucsd.edu}\urladdr{https://klee669.github.io}

\author{Stephen Melczer}
\address{Department of Combinatorics and Optimization, University of Waterloo,	Waterloo, Canada, 200 University Avenue West, Waterloo, ON, CA  N2L 3G1}
\email{smelczer@uwaterloo.ca}\urladdr{https://melczer.ca}

\author{Josip Smolčić}
\address{Department of Combinatorics and Optimization, University of Waterloo,	Waterloo, Canada, 200 University Avenue West, Waterloo, ON, CA  N2L 3G1}
\email{jsmolcic@uwaterloo.ca}\urladdr{https://josip.ca}

\usepackage{latexsym,wasysym, mathrsfs, mathtools,hhline,color, bm}
\usepackage[all, cmtip]{xy}

\usepackage{url}

\usepackage{ textcomp }

\usepackage{graphicx,enumerate}

\usepackage{enumitem}
\usepackage[normalem]{ulem}

\usepackage{framed} 

\usepackage{cite}
\usepackage{amsmath,amssymb,amsfonts}
\usepackage{algorithm,algorithmic}
\usepackage{graphicx}
\usepackage{framed}
\usepackage{textcomp}
\usepackage{xcolor}

\usepackage{multicol}
\usepackage{jlcode}
\usepackage{listings}
\def\jlbasicfont{\ttfamily\selectfont}
\lstdefinestyle{code}{
	language=Julia, 
	showstringspaces=false,
	keywordstyle=\color{blue},
	commentstyle=\color{gray},
	identifierstyle=\color[RGB]{0,102,0},
	columns=fixed,
	keepspaces=false
}
\lstnewenvironment{code}{\lstset{style=code,linewidth=1.01\columnwidth,columns=fullflexible,xleftmargin=0.2cm,framesep=1pt, aboveskip=10pt,belowskip=5pt}}{}

\lstnewenvironment{codeout}{\lstset{style=code,linewidth=1.01\columnwidth,columns=fullflexible,xleftmargin=0.2cm,framesep=1pt, aboveskip=0pt,belowskip=5pt}}{}
\theoremstyle{definition}
	\newtheorem{theorem}{Theorem}[section]
	\newtheorem{lemma}[theorem]{Lemma}

	\newtheorem{definition}[theorem]{Definition}
	\newtheorem{remark}[theorem]{Remark}
	\newtheorem{example}[theorem]{Example}
	
	\def\one{\mathbf{1}} 
	\def\zero{\mathbf{0}} 
	\newcommand{\aaa}{\mathbf{a}}
	\newcommand{\bb}{\mathbf{b}}

	\newcommand{\ii}{\mathbf{i}}

	\newcommand{\pp}{\mathbf{p}}
	\newcommand{\qq}{\mathbf{q}}
	\newcommand{\rr}{\mathbf{r}}
	
	\newcommand{\ww}{\mathbf{w}}
	\newcommand{\xx}{\mathbf{x}}
	\newcommand{\yy}{\mathbf{y}}
	\newcommand{\zz}{\mathbf{z}}

	\newcommand{\HI}{H^{\mfI}}
	\newcommand{\HR}{H^{\mfR}}
	\newcommand{\mfI}{\mathfrak{I}}
	\newcommand{\mfR}{\mathfrak{R}}
	
	\newcommand{\mH}{\mathcal{H}}

	\def\balpha{\boldsymbol{\alpha}}
	
	\def\C{\mathbb{C}} 
	\def\N{\mathbb{N}} 
	\def\Q{\mathbb{Q}} 
	\def\R{\mathbb{R}} 
	\def\Z{\mathbb{Z}} 

\begin{document}
	
	\maketitle
	
	\begin{abstract}
We combine tools from homotopy continuation solvers with the methods of analytic combinatorics in several variables to give the first practical algorithm and implementation for the asymptotics of multivariate rational generating functions not relying on a non-algorithmically checkable `combinatorial' non-negativity assumption. Our homotopy implementation terminates on examples from the literature in three variables, and we additionally describe heuristic methods that terminate and correctly predict asymptotic behaviour in reasonable time on examples in even higher dimension. Our results are implemented in Julia, through the use of the HomotopyContinuation.jl package, and we provide a selection of examples and benchmarks.
	\end{abstract}

\section*{}
Let $(f_n)_{n\in\N}=f_0,f_1,\dots$ be a complex-valued sequence with \emph{generating function} $F(z) = \sum_{n\geq0}f_nz^n$. Although $F$ is a priori only a \emph{formal} power series, in a wide variety of applications (in fact, whenever $f_n$ has at most exponential growth) it represents an analytic function in a neighbourhood of the origin. The field of \emph{analytic combinatorics} creates effective techniques to determine the asymptotic behaviour of $f_n$ through a study of the analytic behaviour of $F(z)$. Most classical methods in analytic combinatorics take as input an algebraic or differential equation satisfied by $F(z)$ and, when successful, return the leading terms in an asymptotic expansion of $f_n$ (see~\cite{FlajoletSedgewick2009} or~\cite[Chapter 2]{Melczer2021}).

More recently, a theory of \emph{analytic combinatorics in several variables (ACSV)}~\cite{Melczer2021,PemantleWilson2013} has been developed to translate the analytic behaviour of a $d$-variate generating function
\[ F(\zz) = \sum_{\ii\in\N^d}f_\ii\zz^\ii := \sum_{\ii\in\N^d}f_{i_1,\dots,i_d}z_1^{i_1}\cdots z_d^{i_d} \]
into asymptotic information about its coefficient sequence $(f_{\ii})_{\ii\in\N^d}$. In this paper we focus on the case of a power series expansion of a multivariate rational function $F(\zz)=G(\zz)/H(\zz)$ and attempt to determine asymptotics of the \emph{$\rr$-diagonal} sequence $(f_{n\rr})_{n\in\N}$ for a fixed \emph{direction vector} $\rr\in\Z_{>0}^d$. The most common situation to arise in practice is the \emph{main diagonal}, when $\rr=\one$.

\begin{remark}
	If $\rr$ has some zero coordinates then we can reduce to the above situation by setting some of the variables equal to zero and working in a lower dimension. For instance, the $(0,r_2,r_3)$-diagonal of any series $F(x,y,z)$ is the $(r_2,r_3)$-diagonal of $F(0,y,z)$. Furthermore, our asymptotic statements continue to hold for directions $\rr\in\Q_{>0}$ if they are interpreted to be valid only when $n\rr\in\N^d$. In fact, the methods of ACSV show that asymptotics of the $\rr$-diagonal usually vary smoothly with $\rr$, allowing one to give a natural interpretation of asymptotics in irrational directions and derive central limit theorems~\cite[Section 5.3.3]{Melczer2021}.
\end{remark}

\begin{remark}
	Because the methods of ACSV hold in any dimension, our requirement that $F(\zz)$ be rational is less restrictive than it may seem. For instance, the $\rr$-diagonal of an algebraic function in $d$ variables can be represented~\cite[Section 3.2.2]{Melczer2021} as the diagonal of a rational function in $2d$ variables (and a `skew-diagonal' of a rational function in $d+1$ variables). The theoretical results discussed here also hold for \emph{meromorphic} functions, when $F(\zz)$ is (locally) the ratio of analytic functions, however our restriction to rational functions allows us to stay in the realm of algebraic quantities and polynomial systems, which we use for our explicit algorithms.
\end{remark}

There are many factors making ACSV more complicated than its univariate counterpart. Although a univariate rational function has a finite number of singularities, meaning one can determine the `asymptotic contribution' of each and simply sum those with the fastest growth, any (non-polynomial) rational function in at least two variables must have an infinite number of singularities. In addition to obscuring which singularities contribute to asymptotics, this also means that the singular set can have \emph{non-trivial geometry}, for instance by self-intersecting. The difficulties that arise mean that unlike the univariate case, which relies on standard complex-analytic results going back hundreds of years, the most advanced ACSV results rely on advanced techniques from areas of mathematics as diverse as complex analysis in several variables, the study of singular integrals, algebraic geometry, differential geometry, and topology.

The starting point of an ACSV analysis expresses the $\rr$-diagonal of $F(\zz)$ as a $d$-dimensional complex integral. In the simplest cases, asymptotic behaviour is determined by the behaviour of $F$ near two types of points: \emph{critical points}, defined by an explicit polynomial system, and \emph{minimal points}, which are singularities that are coordinate-wise closest to the origin. Critical points satisfy a square polynomial system, and generically form a finite set that can be manipulated in a computer algebra system. In contrast, there are always an infinite number of minimal points, which are defined by inequalities involving the moduli of coordinates and are thus trickier and more expensive to manipulate in computations. 

\subsection{Previous Work and Our Contributions}

From the beginning of its modern period in work of Pemantle and Wilson~\cite{PemantleWilson2002}, the goal of ACSV has always been to develop methods explicit enough to be implemented in a computer algebra system. The `surgery' approach of~\cite{PemantleWilson2002}, which applies to generating functions with \emph{smooth} singular sets that form manifolds, essentially computes a residue in one variable to obtain a $(d-1)$-dimensional integral that is approximated using the saddle-point method. Although this surgery method does not require much theory beyond univariate analytic combinatorics, it requires strong conditions on the locations of minimal points that can be computationally expensive to verify. Later techniques, using \emph{cones of hyperbolicity}~\cite{BaryshnikovPemantle2011} and multivariate residue and homology computations~\cite{BaryshnikovMelczerPemantle2022}, rely on more advanced theory but simplify the assumptions that need to be verified for the results to hold. In the simplest cases, which hold for the majority of examples encountered in combinatorial applications, it suffices to determine which of the critical points are minimal and then add explicit asymptotic contributions corresponding to the (finite number of) minimal critical points. The most expensive step in such an analysis is almost always checking minimality.

The first systematic algorithmic study of ACSV methods was conducted by Melczer and Salvy~\cite{MelczerSalvy2021}, who encoded critical points using a symbolic-numeric data structure known as a \emph{Kronecker} or \emph{rational univariate} representation and then reduced checking minimality to rigorously approximating the roots of certain univariate polynomials to sufficiently high accuracy. Those authors created a preliminary implementation of their work, which does not certify numeric computations to provide rigorous proofs and requires \emph{combinatorial} rational functions, in the Maple computer algebra system. A rational function $F(\zz)$ is \emph{combinatorial} if all of its power series coefficients are non-negative: although this condition is satisfied for any multivariate generating function, in many combinatorial examples only one diagonal of $F$ enumerates a combinatorial class and the non-diagonal entries have negative coefficients. It is an open problem, even in the univariate case, whether it is decidable to detect when a rational function is combinatorial (see~\cite{OuaknineWorrell2014} for some open problems in this area). Although Melczer and Salvy~\cite{MelczerSalvy2021} detail a method that, in principle, yields an algorithm for asymptotics that does not require combinatorality, in practice an implementation in Maple would not halt in reasonable time beyond low degree examples in two or three variables.

Instead of continuing with the Kronecker representation approach of Melczer and Salvy, in this paper we exploit homotopy continuation methods to certify minimality of critical points, and ultimately determine asymptotics of $\rr$-diagonals of rational functions. Using the \textsc{HomotopyContinuation.jl} Julia package~\cite{breiding2018homotopycontinuation} for polynomial system solving, we provide the first implementation of ACSV methods under assumptions that often hold in practice. Our implementation is efficient enough to work even without the assumption of combinatorality, although when the user knows a priori that their input rational function is combinatorial then the computation is greatly reduced. In addition, we describe two heuristic methods to classify minimal critical points using numerical approximations that are extremely efficient, and are the only implemented algorithms we currently know of that can aid in the search for minimal points in more than three variables.

\begin{example}
	The main diagonal of the power series expansion of 
	\[ F(x,y,z)=\frac{1}{1-(1+z)(x+y-xy)}\] 
	is related to a result of Apéry~\cite{Apery1983} on the irrationality measure of $\zeta(2)$. After importing our package we define the denominator polynomial in Julia using
	\begin{code}
@polyvar x y z
H = 1-(1+z)*(x+y-x*y)
	\end{code}
	If we know that this power series expansion is combinatorial, then we can get the (truncated for clarity) minimal critical point
	\begin{code}
min_cp = find_min_crits_comb(H)
	\end{code}
	\begin{codeout}
Out: 1-element Vector{Vector{ComplexF64}}:
    [0.38 + e-39im,0.38 + e-38im,0.61 - e-38im]
	\end{codeout}
	and print out the leading asymptotic term of the diagonal with
	\begin{code}
leading_asymptotics(1,H,min_cp)
	\end{code}
	\begin{codeout}
Out: "(0.09+6.2e-39im)^(-n)n^(-1)(0.47-5.7e-40im)"
	\end{codeout}
	It is not obvious from the definition that $F$ is combinatorial. If we don't know our function is combinatorial then we can determine minimality by running \lstinline[style=code]{find_min_crits(H)}, which returns the same point but requires approximately \emph{15 minutes} of computation. If we want to heuristically check for minimal critical points, but don't know that $F$ is combinatorial and don't want to wait for the full algorithm, we can run the algorithms \lstinline[style=code]{find_min_crits(H; approx_crit=true)} or \lstinline[style=code]{find_min_crits(H; monodromy=true)}, described below, which also find the correct point and finish in seconds.
\end{example}

\begin{remark}
	Because we use numeric methods, asymptotic behaviour is returned with numeric approximations of constants. If the user wants to determine the algebraic quantities involved exactly, we recommend solving for the critical points (a relatively cheap operation) symbolically using another computer algebra system like Sage or Maple and then using the results of this package to filter out the minimal ones (the most expensive operation).
\end{remark}

The rest of this paper proceeds as follows. Section~\ref{sec:ACSV} gives a quick recap of the methods of ACSV and the high-level problems that need to be decided to find asymptotics, with a description of numerical algebraic geometry methods for polynomial system solving given in Section~\ref{sec:NAG}. Section~\ref{sec:Julia} uses this background material to detail our \textsc{ACSVHomotopy.jl} Julia package, while Section~\ref{sec:examples} illustrates the package on a wide variety of combinatorial examples, including benchmarks between different algorithms.  Although our algorithms always terminate, due to the nature of homotopy continuation methods they may not always provide a rigorous proof of asymptotics -- Section~\ref{sec:certification} discusses this issue and describes situations in which the algorithms do give rigorous proofs. Finally, Section~\ref{sec:future} concludes with some extensions that we believe should be addressed next.

\section{Smooth ACSV}
\label{sec:ACSV}
From now on, $F(\zz)=G(\zz)/H(\zz)$ denotes a ratio of $d$-variate coprime polynomials $G,H \in \Z[\zz]$ with power series expansion $F(\zz) = \sum_{\ii\in\N^d}f_{\ii}\zz^{\ii}$ converging around the origin, and $\rr \in \Z_{>0}^d$ is a fixed direction vector.

\begin{definition}[minimal critical points] 
	A point $\ww \in \C_*^d$ is a \emph{(simple) smooth critical point} of $F$ if $(\nabla H)(\ww)\neq\zero$ and
	\begin{equation}
		\left\{\begin{array}{l}
			H(\ww)=0\\
			r_kz_1H_{z_1}(\ww) - r_1z_kH_{z_k}(\ww) = 0 \;~(2\leq k\leq d).
		\end{array}\right. 	\label{eq:CP}
	\end{equation}
	We call $\ww \in \C_*^d$ a \emph{minimal point} if $H(\ww)=0$ and there does not exist $\yy\in\C^d$ such that $H(\yy)=0$ and $|y_j|<|w_j|$ for all $j=1,\dots,d$.
\end{definition}

\begin{remark}
	If $(\nabla H)(\ww)=\zero$ then~\eqref{eq:CP} is trivially satisfied. If the gradient vanishes because $H$ has a higher-order pole (for instance, if $H=P^2$ for some polynomial $P$) then our analysis of minimal critical points can be performed on the square-free part of $H$ (the product of its irreducible factors) to obtain an asymptotic expansion of $f_{n\rr}$ with minor modifications. On the other hand, if the gradient vanishes because the zero set of $H$ self-intersects then more advanced techniques are required~\cite[Part III]{Melczer2021}.
\end{remark}

We will be able to determine asymptotics in the presence of smooth minimal critical points, assuming a nondegeneracy condition on the zero set of $H$.

\begin{definition}[phase Hessian matrix] 
	If $\ww$ is a smooth critical point then the \emph{phase Hessian matrix} $\mH$ at $\ww$ is the $(d-1)\times(d-1)$ matrix defined by
	\[
	\mH_{i,j} = 
	\begin{cases}
		V_iV_j + U_{i,j} - V_jU_{i,d} - V_iU_{j,d} + V_iV_jU_{d,d} &: i \neq j \\[+3mm]
		V_i + V_i^2 + U_{i,i} - 2V_iU_{i,d} + V_i^2U_{d,d} &: i=j
	\end{cases}
	\]
	where
	\[ U_{i,j} = \frac{w_iw_j H_{z_iz_j}(\ww)}{w_dH_{z_d}(\ww)} \qquad \text{and} \qquad V_i = \frac{r_i}{r_d}.\]
\end{definition}

\begin{theorem}[{Melczer~\cite[Theorem 5.1]{Melczer2021}}]
	\label{thm:ACSV}
	Suppose that the system of polynomial equations~\eqref{eq:CP} admits a finite number of solutions, exactly one of which, $\ww\in\C_*^d$, is minimal. Suppose further that $H_{z_d}(\ww)\neq0$, that the phase Hessian matrix $\mH$ at $\ww$ has non-zero determinant, and that $G(\ww) \neq 0$. Then, as $n\rightarrow\infty$,
	
	{\small
		\[
		f_{n\rr} = \ww^{-n\rr} n^{(1-d)/2} \frac{(2\pi r_d)^{(1-d)/2}}{\sqrt{\det(\mH)}} \frac{-G(\ww)}{w_d\, H_{z_d}(\ww)}\left(1 + O\left(\frac{1}{n}\right)\right).
		\]
	}
	
	When the zero set of $H$ contains a finite number of points with the same coordinate-wise modulus as $\ww$, all of which satisfy the same conditions as $\ww$, then an asymptotic expansion of $f_{n\rr}$ is obtained by summing the right hand side of this expansion at each point.
\end{theorem}

\begin{remark}
	The condition that $G(\ww)\neq0$ means that the leading asymptotic term in Theorem~\ref{thm:ACSV} doesn't vanish. When $G(\ww)=0$ asymptotics can usually still be determined by computing higher-order terms using (increasingly complicated) explicit formulas.
\end{remark}

\subsection{Minimality Tests}

The hardest work in applying Theorem~\ref{thm:ACSV} is computing the critical points, defined implicitly by~\eqref{eq:CP}, and determining which, if any, are minimal. 

\textbf{(Combinatorial Case)}
Recall that a function is called \emph{combinatorial} if its power series expansion contains only a finite number of negative coefficients. When $F$ is combinatorial there is a simple test for minimal critical points.

\begin{lemma}[Melczer and Salvy~\cite{MelczerSalvy2021}]
	\label{lem:combmin}
	Suppose $F$ has only a finite number of negative power series coefficients $f_\ii$. If $\yy\in\C_*^d$ is a minimal critical point then so is $(|y_1|,\dots,|y_d|)$. Furthermore, $\ww \in \R_{>0}^d$ is a minimal critical point if and only if the system
	\begin{equation}
		\begin{split}
			H(\zz) = H(tz_1,\dots,tz_d) &= 0 \\
			z_1H_{z_1}(\zz)-r_1\lambda = \cdots = z_dH_{z_d}(\zz)-r_d\lambda &= 0
		\end{split}
		\label{eq:extendedSys}
	\end{equation}
	has a solution $(\zz,\lambda,t) \in \R^{d+2}$ with $\zz=\ww$ and $t=1$ \emph{and} no solution with $\zz=\ww$ and $0 < t < 1$.
\end{lemma}

In the combinatorial case we firstly use Lemma~\ref{lem:combmin} to characterize the minimal critical points with positive coordinates by studying the solutions to~\eqref{eq:extendedSys}. From them, find the solutions to~\eqref{eq:CP} with the same coordinate-wise modulus. The following algorithm summarizes this approach.

\begin{algorithm}	[H]
	\caption{Minimal Critical Points in the Combinatorial Case}
    \label{algo:comb}
	\begin{enumerate}
		
		\item Determine the set $S$ of zeros of the polynomial system (\ref{eq:extendedSys}) in the variables $\mathbf{z},\lambda,t$. If $S$ is not finite, FAIL.
		\item Find $\boldsymbol{\zeta}\in \mathbb{R}^d_{>0}$ such that there exists $(\boldsymbol{\zeta},\lambda,t)\in S$ and for all such triples, $t\not\in (0,1)$. If the number of such $\boldsymbol{\zeta}$'s is not exactly $1$ or if there are such points with $\lambda=0$, FAIL.
		\item Identify $\boldsymbol{\zeta}$ among the elements of the set $\mathcal{C}$ of zeros to (\ref{eq:CP}).
		\item Return $$\{\mathbf{z}\in \mathbb{C}^d\mid \exists (\mathbf{z},\lambda)\in \mathcal{C},~|z_1|=|\zeta_1|,\cdots,|z_d|=|\zeta_d|\}.$$
	\end{enumerate}
\end{algorithm}

\textbf{(General Case)}
If $F$ is not combinatorial, or if we don't know a priori that $F$ is combinatorial, then it is no longer sufficient to consider only the critical points with positive real coordinates to check minimality. In order to express the moduli of coordinates as algebraic equations, we write $H(\xx+i\yy) = \HR(\xx,\yy) + i\HI(\xx,\yy)$ for real variables $\xx,\yy\in\R^d$ and polynomials $\HR,\HI \in \R[\xx,\yy]$. Translating the smooth critical point equations~\eqref{eq:CP} into these new coordinates gives that $\zz=\aaa+i\bb$ with $\aaa,\bb\in\R^d$ is critical if and only if
\begin{align}
	\HR(\aaa,\bb) = \HI(\aaa,\bb) &= 0  \label{eq:GenSys1} \\
	a_j \HR_{x_j}(\aaa,\bb) + b_j\HR_{y_j}(\aaa,\bb) - r_j\lambda_R&=0  \label{eq:GenSys2} \\
	a_j \HI_{x_j}(\aaa,\bb) + b_j\HI_{y_j}(\aaa,\bb) - r_j\lambda_I&=0  \label{eq:GenSys3} 
\end{align}
for some $\lambda_R,\lambda_I \in R$, where $1 \leq j \leq d$ in each equation. To test minimality of these critical points we add the equations
\begin{align}
	\HR(\xx,\yy) = \HI(\xx,\yy) &= 0 \label{eq:GenSys4} \\
	x_j^2 + y_j^2 - t(a_j^2+b_j^2) &= 0 \label{eq:GenSys5}
\end{align}
for $1 \leq j \leq d$, and verify there is no real solution to~\eqref{eq:GenSys1}-\eqref{eq:GenSys5} with $0 < t < 1$. Generically~\eqref{eq:GenSys1}-\eqref{eq:GenSys5} have a finite set of \emph{real} solutions, corresponding to the generically finite number of critical points of $F$, but because this system contains $3d+4$ equations in $4d+3$ variables it will never have a (non-zero) finite number of solutions over the complex numbers. By considering critical values of the projection map onto the $t$ coordinate, Melczer and Salvy~\cite{MelczerSalvy2021} proved that minimality can be tested by adding the additional equations
\[ (\nu_1y_j-\nu_2 x_j)\HR_{x_j}(\xx,\yy) - (\nu_1x_j+\nu_2 y_j)\HR_{y_j}(\xx,\yy) =0 \]
for $1 \leq j \leq d$. When $\nu_1 \neq 0$ then we can scale by $\nu_1$ and introduce the equations
\begin{equation}
	(y_j-\nu x_j)\HR_{x_j}(\xx,\yy) - (x_j+\nu y_j)\HR_{y_j}(\xx,\yy) =0
	\label{eq:GenSys6}
\end{equation}
to~\eqref{eq:GenSys1}-\eqref{eq:GenSys5}, resulting in a square system with $4d+4$ variables and equations. The case when $\nu_1=0$ is dealt with separately by adding the equations 
\begin{equation*}
	\hspace{0.65in}
	-x_j\HR_{x_j}(\xx,\yy) - y_j\HR_{y_j}(\xx,\yy) = 0
	\hspace{0.65in} 
	(\ref{eq:GenSys6}')
\end{equation*}
for $1 \leq j \leq d$. We determine the minimal critical points by finding $\pp + i\qq$ such that equations \eqref{eq:GenSys1}-\eqref{eq:GenSys3} have a real solution with $(\aaa,\bb)=(\pp,\qq)$ but neither \eqref{eq:GenSys1}-\eqref{eq:GenSys6} nor \eqref{eq:GenSys1}-(\ref{eq:GenSys6}') have a real solution with $(\aaa,\bb)=(\pp,\qq)$ and $0<t<1$.
This process provides the following algorithm.
\begin{algorithm}[H]
	Minimal Critical Points in the Non-Combinatorial Case
	\begin{enumerate}
		
		\item Determine the set $S$ of zeros of the polynomial system (\ref{eq:GenSys1})-(\ref{eq:GenSys6}) in the variables $\mathbf{a},\mathbf{b},\mathbf{x},\mathbf{y},\lambda_R,\lambda_I,\nu,t$. If $S$ is not finite, FAIL.
		\item Construct a set $\mathcal{U}$ of minimal critical points $\mathbf{a}+i\mathbf{b}\in \mathbb{C}^d$ such that there exists $(\mathbf{a},\mathbf{b},\mathbf{x},\mathbf{y},\lambda_R,\lambda_I,\nu,t)\in S\cap \mathbb{R}^{4d+4}$ and for all such tuples, $t\not\in (0,1)$. If either $\mathcal{U}$ is empty or one of its elements has $\lambda_I=\lambda_R=0$, or if the elements of $\mathcal{U}$ do not all belong to the same torus, FAIL.
		\item Identify the elements of $\mathcal{U}$ within the set $\mathcal{C}$ of zeros to (\ref{eq:CP}) and return them.
		\item Do the same for the polynomial system (\ref{eq:GenSys1})-(\ref{eq:GenSys6}') in the variables $\mathbf{a},\mathbf{b},\mathbf{x},\mathbf{y},\lambda_R,\lambda_I,t$.
	\end{enumerate}
\end{algorithm}

Unfortunately, to moving $4d+4$ variables makes verifying minimality much less practical than the combinatorial case. In essence, Lemma~\ref{lem:combmin} states that to prove minimality in the combinatorial case it is sufficient to consider specific line segments in $\R^d$, while to prove minimality in the general case one must consider a much larger set of points in $\C^d$ whose coordinate-wise moduli lie on specific line segments in $\R^d$.

\begin{remark}
	Melczer and Salvy~\cite{MelczerSalvy2021} incorrectly state that $\nu_1$ and $\nu_2$ must both be non-zero: at least one is non-zero at the solutions of interest, but the other may vanish. This is why we introduce~(\ref{eq:GenSys6}').
	Melczer and Salvy~\cite{MelczerSalvy2021} also require an extra condition that a certain Jacobian matrix is non-singular, however this is mainly required for their complexity analysis. If this condition fails then the system~\eqref{eq:GenSys1}-\eqref{eq:GenSys6} can have extra solutions that are irrelevant to detecting minimality, but the presence of such solutions does not affect correctness of the minimality test.
\end{remark}

\section{Numerical Algebraic Geometry}
\label{sec:NAG}
Having reduced the ACSV analysis to questions about polynomial systems, we now recall some methods in computational algebraic geometry for the study of such systems. Although the theory of Gr\"obner bases is, by now, the basis of much work in this area, more recently \emph{numerical algebraic geometry} has emerged as a practical alternative. In this section, we discuss several topics in numerical algebraic geometry that will be used for our techniques.

\subsection{Homotopy Continuation}
\label{sec:homotopy}
\emph{Homotopy continuation} is one method to find numerical approximations of solutions to an $n\times n$ square system $\mathcal{F}=(f_1,\dots, f_n)$ of polynomial equations with $n$ variables. From the system $\mathcal{F}$ we construct an $n\times n$ polynomial system $\mathcal{G}$ whose solutions are known a priori. The system $\mathcal{G}$ is called a \emph{start system} and the system $\mathcal{F}$ is called the \emph{target system}: connecting $\mathcal{F}$ and $\mathcal{G}$ using a homotopy $\mathcal{H}(x,t)$ such that $\mathcal{H}(x,0)=\mathcal{G}$ and $\mathcal{H}(x,1)=\mathcal{F}$, we obtain solutions of $\mathcal{F}$ by tracking homotopy paths from $t=0$ to $t=1$.
To track the homotopy paths, a numerical ordinary differential equation solving technique called the \emph{Davidenko equation} and Newton iteration are used. These tracking techniques are typically referred to as \emph{predictor-corrector} methods. For details, see \cite[Chapter 2]{SommeseWampler2005}. Homotopy continuation is implemented in \textsc{Bertini} \cite{BHSW06}, \textsc{HomotopyContinuation.jl} \cite{breiding2018homotopycontinuation}, and \textsc{NAG4M2} \cite{leykin2011numerical}.

\subsection{Polyhedral Homotopy Continuation}
\label{sec:phc}

The complexity of solving a polynomial system using homotopy continuation is determined by the number of homotopy paths to track. It is thus important to track a number of paths that are at least as large as the number of solutions of the system (so that all solutions can be found) but is not too much larger (to save computation). For polytopes $Q_1,\dots, Q_n$ the Euclidean volume $\text{Vol}(a_1Q_1+\cdots +a_nQ_n)$ of the \emph{Minkowski sum} $a_1Q_1+\cdots +a_nQ_n$ is a homogeneous polynomial in the $n$ variables $a_1,\dots, a_n$, whose coefficient of $a_1a_2\cdots a_n$ is the \emph{mixed volume} $\text{MVol}(Q_1,\dots, Q_n)$ of $Q_1,\dots, Q_n$. 


\begin{theorem}[{Bernstein's theorem~\cite[Theorem A]{bernshtein1975number}}]
	Let $\mathcal{F}$ be a system of polynomials $f_1,\dots, f_n$ in $\C[x_1,\dots, x_n]$. The number of isolated solutions of $\mathcal{F}$ over $\C_*^n$ is at most $\text{MVol}(Q_{f_1},\dots, Q_{f_n})$, where $Q_{f_i}$ is the Newton polytope of $f_i$. Furthermore, for polynomials $f_1,\dots, f_n$ with generic coefficients, the number of solutions for $\mathcal{F}$ over the torus is exactly $\text{MVol}(Q_{f_1},\dots, Q_{f_n})$.
\end{theorem}

The \emph{polyhedral homotopy continuation} method established by Huber and Sturmfels \cite{huber1995polyhedral} is one common way to construct a start system whose solutions form a set with the size of the mixed volume of a system. Consider a polynomial 
\[f(\xx) = \sum\limits_{\aaa\in A}c_\aaa\xx^\aaa\in \C[x_1,\dots, x_n]\]
where $A$ is a collection of integer lattice points. Multiplying each monomial $\xx^\aaa$ of $f$ by some term $t^{w(\aaa)}$ for a \emph{lifting function} $w:A\rightarrow \Z$, we obtain the \emph{lifted polynomial}
\[\overline{f}(\xx,t)=\sum\limits_{\aaa\in A}c_\aaa\xx^\aaa t^{w(\aaa)}.\]
Suppose that a target system $\mathcal{F}$ consists of polynomials $f_1,\dots, f_n$ supported on $A_{f_1},\dots, A_{f_n}$, respectively. Lifting all polynomials $f_1,\dots, f_n$ in $\mathcal{F}$ gives a lifted system $\overline{\mathcal{F}}(\xx,t)$ satisfying $\overline{\mathcal{F}}(\xx,1)=\mathcal{F}$. The solutions of $\overline{\mathcal{F}}$ can be expressed by Puiseux series $\xx(t)=(x_1(t),\dots, x_n(t))$ where 
\[x_i(t)=t^{\alpha_i}y_i+ \text{ higher order terms}\]
for some $\alpha_i\in\Q$ and nonzero constant $y_i$, and substituting $\xx(t)$ back into our polynomials gives
\[\overline{f}_j(\xx(t),t) =\sum\limits_{\aaa\in A_{f_j}}c_\aaa \yy^\aaa t^{\langle \aaa,\balpha\rangle+w(\aaa)}+\text{ higher order terms}.\] 
For a suitable choice of $w$, the constants $\yy$ and exponents $\aaa$ can be computed at each branch of $\overline{\mathcal{F}}$, ultimately describing a start system $\mathcal{G}=\overline{\mathcal{F}}(\xx,0)$ with the right number of solutions. The polyhedral homotopy continuation is implemented in \textsc{HOM4PS2} \cite{lee2008hom4ps}, \textsc{HomotopyContinuation.jl} \cite{breiding2018homotopycontinuation}, and \textsc{PHCpack} \cite{verschelde1999algorithm}.

\subsection{Monodromy}

As seen in Section~\ref{sec:ACSV}, we typically have some solutions of a polynomial system representing critical points and want to determine additional solutions to rule out those that are non-minimal. This `bootstrapping' can be accomplished by monodromy.

For $m,n\in \mathbb{N}$, consider the complex linear space of $n\times n$ square systems $\mathcal{F}_p=(f_p^1,\dots, f_p^n)$ depending on some coefficient parameters $p\in \C^m$, where the monomial support for each polynomial $f_p^i$ is fixed. If we consider an affine linear map $\varphi:p\mapsto F_p$ for $p\in\mathbb{C}^m$ then we can write $\varphi(\mathbb{C}^m)=B$, where $B$ is a parametrized linear variety of systems, and we define the \emph{solution variety} $V=\{(F_p,x)\in B\times \mathbb{C}^n\mid F_p(x)=0\}$ and \emph{projection map} $\pi:V\rightarrow B$. 

Assume that the fiber $\pi^{-1}(F_p)$ only has finitely many points for a generic choice of $p$. The set $D$ of systems in $B$ with non-generic fiber is called the \emph{branch locus} of $\pi$. Each element in the \emph{fundamental group} $\pi_1(B\setminus D)$ of loops in $B\setminus D$ modulo homotopy equivalence induces a permutation on the fiber $\pi^{-1}(F_p)$, which is called a \emph{monodromy action}. To find all solutions of a system $F_p\in \pi(V)$ with generic $p$, one can first find a seed solution $(p_0,x_0)\in V$ and numerically compute the monodromy action to find all solutions of $F_p$. When the solution variety $V$ is irreducible then the monodromy action is transitive. This method for finding solutions of polynomial systems is studied and implemented in \cite{breiding2018homotopycontinuation,duff2019solving}.

\subsection{Certification}

By construction, numerical methods return approximations, so some kind of certification is necessary for rigorous results. Specifically, a user needs a certificate that an approximation obtained by the homotopy method is properly approximating a solution of a system. A numerical approximation is called \emph{certified} if it can be refined to an actual solution of the system to an arbitrary precision by applying iterative operators (such as Newton iteration). Software providing such certification includes \textsc{alphaCertified} \cite{hauenstein2011alphacertified}, the function \textsc{certify} implemented in \textsc{HomotopyContinuation.jl} \cite{breiding2020certifying} and \textsc{NumericalCertification} \cite{https://doi.org/10.48550/arxiv.2208.01784}. In our implementation, we use the function \textsc{certify} in \textsc{HomotopyContinuation.jl} exploiting the Krawczyk's method via interval arithmetic \cite[Chapter 8]{moore2009introduction}.

\section{The ACSVHomotopy Package}
\label{sec:Julia}

We now combine the theory of ACSV presented in Section~\ref{sec:ACSV} with the techniques described in Section~\ref{sec:NAG} to create effective and practical algorithms for the asymptotics of multivariate rational functions. Our algorithms are implemented in the \textsc{Julia} package \textsc{ACSVHomotopy.jl}, using the \textsc{HomotopyContinuation.jl} package for our homotopy and monodromy computations. 

The package is available at 
\begin{center}
	{github.com/ACSVMath/ACSVHomotopy}
\end{center}
and our example worksheet can be viewed at 

\begin{center}
	{github.com/ACSVMath/ACSVHomotopy/blob/main/ExampleWorksheet.ipynb}
\end{center}

\subsection{Combinatorial Case}

For the combinatorial case we first compute the distinct solutions to~\eqref{eq:CP} using a polyhedral homotopy with certification by Krawczyk's method. We then solve and certify~(\ref{eq:extendedSys}) with the added equation $(1-t)\mu-1=0$ to eliminate all solutions with $t=1$ (there are never any solutions with $t=0$ as this would imply $H(\zero)=0$, contradicting $F$ having a power series expansion). Since we no longer have solutions where $t=1$, by refining the solutions to sufficient precision we can determine the solutions with positive real coordinates where $0<t<1$, match the projection onto the $\zz$ variables of each to a distinct solution of~\eqref{eq:CP}, and thus rule out all non-minimal critical points with positive coordinates. We then find all critical points with the same coordinate-wise moduli and return that set. 

\begin{example}
	As a simple example, we can find the minimal critical point
	\begin{code}
@polyvar x y
find_min_crits_comb(1-x-y)
	\end{code}
	\begin{codeout}
Out:1-element Vector{Vector{ComplexF64}}:
    [0.5 + 0.0im, 0.5 + 0.0im]
	\end{codeout}
	controlling asymptotics for the central binomial coefficient
	$\binom{2n}{n}$ which forms the main diagonal sequence of 
	\[ F(x,y) = \frac{1}{1-x-y}.\]
	Similarly, we can compute the approximations
	\begin{code}
@polyvar x y z
find_min_crits_comb(1-z*(x^2*y+y+x*y^2+x))
	\end{code}
	\begin{codeout}
Out:2-element Vector{Vector{ComplexF64}}:
    [1.0 + e-35im, 1.0 - e-35im, 0.25 - e-37im]
    [-1.0 - e-36im, -1.0 + e-36im, -0.25 - e-36im]
	\end{codeout}
	for the two minimal critical points $\pm(1,1,1/4)$ determining
	asymptotics for the main diagonal of 
	\[F(x,y,z) = \frac{(1+x)(1+y)}{1-zxy(x+1/x+y+1/y)}.\] 
	This diagonal enumerates walks on the cardinal directions $\{N,S,E,W\}=\{(\pm1,0),(0,\pm1)\}$ that start at the origin and stay in $\N^2$.
\end{example}

As in the work of Melczer and Salvy, the most expensive operation occurs when trying to group roots with the same coordinate-wise modulus as a known minimal critical point (step~(4) in Algorithm~\ref{algo:comb}). When using a symbolic-numeric method it is possible to compute minimal polynomials for the values of the coordinates and use this to identify points with the same coordinate-modulus by computing numerically to $O(h\delta^{3d})$ bits of precision, where $h$ is a bound on the bitsize of the coefficients of the denominator $H$ and $\delta$ is the degree of $H$ (see~\cite[Corollary 54]{MelczerSalvy2021}). Combining past bounds in the literature~\cite{Safey-El-DinSchost2018,DubickasSha2015} allows us to identify an explicit precision such that if two solutions have coordinate-wise moduli agreeing to this precision then their coordinate-wise moduli are exactly equal. Our bound is also of $O(h\delta^{3d})$ bits, however the constant in front is worse than that when using minimal polynomials. After a minimal critical point is identified, our code continually refines precision using Newton iteration until any points with the same coordinate-wise moduli are found. In practice, any points with different coordinate-wise moduli are identified at precision much lower than the worst case bound, but we must compute up to the bound when there are points with the same coordinate-wise moduli. Unfortunately, the extreme precision involved with a large number of variables means that we cannot always rigorously check to the required accuracy, and the code returns a warning to the user with its output when this occurs.

\subsection{General Case}
In general we must consider the extended systems~\eqref{eq:GenSys1}-\eqref{eq:GenSys6} and~\eqref{eq:GenSys1}-(\ref{eq:GenSys6}'), which essentially doubles the number of variables under consideration. Mirroring the combinatorial case, we can solve~\eqref{eq:GenSys1}-\eqref{eq:GenSys6} and~\eqref{eq:GenSys1}-(\ref{eq:GenSys6}') using a polyhedral homotopy, certify the results using Krawczyk's method, and refine to a sufficient precision to determine when $0 < t < 1$ to rule out non-minimal points. 

\begin{remark}
	The system~\eqref{eq:GenSys1}-(\ref{eq:GenSys6}') is over determined, with $4d+4$ equations and $4d+3$ variables. In order to use \textsc{HomotopyContinuation.jl} we drop one of the equations in~(\ref{eq:GenSys6}') to obtain a square system: this can introduce additional solutions which are irrelevant to determining minimality, but does not affect the correctness of our test for minimality.
\end{remark}

\begin{example}
	Straub and Zudilin~\cite{StraubZudilin2015}, following Gillis, Reznick, and Zeilberger~\cite{GillisReznickZeilberger1983}, study families of rational functions connected to special function theory. For instance, in three dimensions they study the constants $c$ for which 
	\[ F_c(x,y,z) = \frac{1}{1-(x+y+z)+cxyz}\] 
	has non-negative power series coefficients on its main diagonal (which turns out to imply non-negativity of all power series coefficients). Running the code (for $c=5$)
	\begin{code}
@polyvar x y z
find_min_crits(1 - (x+y+z) + 5*x*y*z)
	\end{code}
	\begin{codeout}
Out:2-element Vector{Vector{ComplexF64}}:
    [0.45 - 0.12im, 0.45 - 0.12im, 0.45 - 0.12im]
    [0.45 + 0.12im, 0.45 + 0.12im, 0.45 + 0.12im]
	\end{codeout}
	gives the minimal critical points controlling asymptotics of the main diagonal. Since this is a complex conjugate pair, the resulting asymptotic expansion implies that $F_c$ has an infinite number of negative coefficients on its main diagonal when $c=5$ (in fact, it has an infinite number of negative coefficients whenever $c>4$).
\end{example}

\subsection{Faster Heuristics}

As seen in the examples of Section \ref{sec:examples}, the high number of variables in the extended system even in low dimensional examples means it does not terminate within reasonable time for polynomials with four or more variables. In order to speed up our solvers, we can numerically approximate the distinct solutions to the small system~\eqref{eq:GenSys1}-\eqref{eq:GenSys3} and then substitute each of these solutions as parameters into the extended equations~\eqref{eq:GenSys4}-\eqref{eq:GenSys6} and~\eqref{eq:GenSys4}-(\ref{eq:GenSys6}'). In the implementation, it is done by running the function \lstinline[style=code]{find_min_crits} with the flag \lstinline[style=code]{approx_crit = true}.

\begin{remark}
	This approach approximates the solutions to the extended system~\eqref{eq:GenSys1}-\eqref{eq:GenSys6} if the solutions vary smoothly with $\aaa$ and $\bb$, which happens whenever the Jacobian of~\eqref{eq:GenSys4}-\eqref{eq:GenSys6} with respect to $\aaa$ and $\bb$ is full rank at all values of $\aaa$ and $\bb$ solving \eqref{eq:GenSys1}-\eqref{eq:GenSys3}. Unfortunately, verifying this condition is usually about as costly as solving the extended system, so we do not do this in our computations and refer to this method only as an efficient heuristic that correctly identifies minimal critical points in a large variety of cases.
\end{remark}

\begin{example}
	To stress-test our algorithms we generate a random polynomial $p(x,y,z)$ with six terms in four variables having coefficients in $\{1,\dots,100\}$ and then set $H(x,y,z)=1-p(x,y,z)$. Running
	\begin{code}
@polyvar x y z
H=1-(72*x^3*z+97*y*z^3+53*x*z^2+47*x*y+39*z^2+71*x)
find_min_crits(H; approx_crit = true)
	\end{code}
	\begin{codeout}
Out:1-element Vector{Vector{ComplexF64}}:
    [0.001+5.5e-40im, 6.2-7.5e-37im, 0.06+0.0im]
	\end{codeout}
	returns the unique minimal critical point in about three minutes. This example does not terminate without the \lstinline[style=code]{approx_crit = true} flag.
\end{example}

It is also possible to use the approximations to the critical points as a start system to solve~\eqref{eq:GenSys4}-\eqref{eq:GenSys6} using the monodromy method. More precisely, for any $(\mathbf{a},\mathbf{b})$ solving \eqref{eq:GenSys1}-\eqref{eq:GenSys3} we set $(\xx,\yy)=(\mathbf{a},\mathbf{b})$ and $t=1$ in~\eqref{eq:GenSys4}-\eqref{eq:GenSys6} and then compute a corresponding start value of $\nu$ by computing the left kernel of the Jacobian matrix of \eqref{eq:GenSys4}-\eqref{eq:GenSys5} with respect to variables $\mathbf{x},\mathbf{y}$ and $t$. From a given parameter value $(\mathbf{a},\mathbf{b})$ and the initial solution $(\mathbf{x},\mathbf{y},\nu,t)$, we collect real solutions from the monodromy method and check if $t\in (0,1)$: if it is then we remove the parameter value $(\mathbf{a},\mathbf{b})$  as it is non-minimal. Interestingly, it appears that monodromy cannot detect solutions where $\nu=0$ when starting with a non-zero value of $\nu$, and vice-versa, suggesting that solution variety is the union of components corresponding to these cases. We thus repeat this process separately for the cases where $\nu=0$ and when~(\ref{eq:GenSys6}') replaces~\eqref{eq:GenSys6}. Finally, we return the values of $(\mathbf{a},\mathbf{b})$ that are not disregarded.

\begin{example}
	Melczer and Salvy~\cite{MelczerSalvy2021} introduce the rational function
	\[ F(x,y) = \frac{1}{(1-x-y)(20-x-40y)-1}  \]
	because it has two critical points with positive coordinates, one of which is smaller in the first coordinate and the other of which is smaller in the second coordinate (so it is not clear which, if any, should be minimal). Running
	\begin{code}
@polyvar x y
H = (1-x-y)*(20-x-40*y)-1 
find_min_crits(H; monodromy=true)
	\end{code}
	\begin{codeout}
Out:1-element Vector{Vector{ComplexF64}}:
    [0.54 - 9.18e-41im, 0.31 + 1.83e-40im]
	\end{codeout}
	returns the correct minimal critical point.
\end{example}

We believe that further study of the geometric properties of the extended system~\eqref{eq:GenSys4}-\eqref{eq:GenSys6} could help make this monodromy approach a powerful tool for ACSV analysis.

\section{Examples and Benchmarks}
\label{sec:examples}

Tables~\ref{table1} and~\ref{table2} list benchmarks of our implementation against a selection of combinatorial and algebraic examples, executed on a {Macbook pro, 2 GHz Quad-Core Intel Core i5, 16 GB RAM}. The package supports arbitrary $\rr$-diagonal sequences, but examples in this section were done with $\rr=\one$. See our supplementary notebook for the full details on the rational functions involved.

\begin{remark}
	The \textsc{HomotopyContinuation.jl} package converts input polynomials into compiled straight-line programs for fast evaluation. In order to better see the differences between examples as they grow in degree and dimension, we have removed compilation time from our benchmarks (compilation time takes the majority of the runtime for small examples but is a small part of larger examples). This results in several seconds on small examples, and up to tens of seconds on larger examples, that are not included in the reported timings. In particular, the (non-certified) package of Melczer and Salvy beats our package in the combinatorial case on most examples in Table~\ref{table1} when compilation time is added (except for the two high degree examples where the Maple package takes much longer).
\end{remark}

\begin{table}
	\centering
	\begin{tabular}{c|c|c}
		Example & Comb. & Maple Comb.  \\
		\hline
		$1-x-y$  & 0.0052 & 0.143  \\
		Two positive CPs & 0.029 & 0.292  \\
		square-root & 0.01  & 0.06 \\
		Apéry $\zeta(2)$ & 0.025 & 0.06 \\
		Apéry $\zeta(3)$ & 0.7  & 0.3 \\ 
		random poly & 0.9  & 840 \\
		2D Walk & 0.03 & 0.06 \\
		3D Walk & 0.08 & 2.7 \\
		$1-x-y^2-w^3-z^4$ & 0.06 & 509
	\end{tabular}
	\vspace{0.1in}
	\caption{Time, in seconds, of running our Julia implementations in the combinatorial case,
		compared to the Kronecker representation approach of Melczer and Salvy. The time to compile Julia functions is not included. 
	}
	\label{table1}
\end{table}

\begin{table}
	\centering
	\begin{tabular}{c|c|c|c}
		Example & HSolve & HSolve Approx & Monodromy  \\
		\hline
		$1-x-y$ & 0.04 & 0.02 & 2.3 \\
		Two positive CPs & 4.1 & 0.33 & 2.84 \\
		Apéry $\zeta(2)$ & 670 & 3.8 & 8.5  \\
		square-root & 29.5 & 0.72 & 14.9 \\
		random poly & INC & 189.4 & 583.1  \\
		2D Walk & INC & 15.3 & 31.9 \\
		GRZ & 236 & 3.6 & 3.8
	\end{tabular}
	\vspace{0.1in}
	\caption{Time, in seconds, of running our Julia implementations that do not assume combinatorality. The first column is the time to solve the extended critical point systems, the second column is the time to solve the smaller systems after solving the critical point system separately, and the final column is the time to run the monodromy method. INC indicates the code did not complete after running for an hour.}
	\label{table2}
\end{table}

\section{Rigor of Results}
\label{sec:certification}

Because we certify our solutions, we never attempt to approximate a point that is not actually a solution of the polynomial systems under consideration. However, by the way they are designed, it is possible for homotopy computations to \emph{miss} solutions, which could result in a point being deemed minimal when it is not. There are some exceptions: when the number of solutions found matches the upper bound on the number of solutions given by the mixed volume, for instance, then we can be sure we have found all solutions. Tables~\ref{tab:mv1} and~\ref{tab:mv2} show a comparison between the mixed volumes of several systems studied here, compared to the actual number of solutions found. It can be observed that we often reach the upper bound in the combinatorial case, but this usually does not happen in the non-combinatorial case.

\begin{table}
	\centering
	\begin{tabular}{c|c|c}
		Example & Mixed volume & $\#$ Solutions found  \\
		\hline
		$1-x-y$  &1 & 1 \\
		$1 - xy - xy^2 - 2x^2y
		$ & 9 & 9 \\
		$1 - x - y^2 - w^3 - z^4$  &  96 & 96 
	\end{tabular}
	\vspace{0.1in}
	\caption{Mixed volume for the system (\ref{eq:extendedSys}) and the actual number of solutions found for several combinatorial examples}
	\label{tab:mv1}
\end{table}

\begin{table}
	\centering
	\begin{tabular}{c|c|c|c|c}
		Example & (\ref{eq:GenSys1})-(\ref{eq:GenSys6}) & $\#$ Sols  & (\ref{eq:GenSys1})-(\ref{eq:GenSys6}') & $\#$ Sols \\
		\hline
		$1-x-y$  &4 & 1 & 2 & 0 \\
		$1 - xy - xy^2 - 2x^2y
		$ & 3276 & 99 & 1638 & 126 \\
		$1 - (x+y+z) + \frac{81}{8}xyz
		$  &  13068 & 162 & 4356 & 216 \\
		$1 - x - y^2 - w^3 - z^4$  & FAIL &N/A &  442368 & 442368
	\end{tabular}
	\vspace{0.1in}
	\caption{Mixed volumes for the systems (\ref{eq:GenSys1})-(\ref{eq:GenSys6}) and (\ref{eq:GenSys1})-(\ref{eq:GenSys6}') and the number of solutions found for several examples. FAIL means that the code did not compute the mixed volume due to an out of memory error.}
	\label{tab:mv2}
\end{table}

We can also conclude we know minimality rigorously when there is another way to determine that a minimal critical point must exist, and all but one point is ruled out by our algorithms. For instance, in the combinatorial case it can be shown that any polynomial whose support contains the terms $1,z_1,z_2,\dots,z_d$ must have at least one minimal critical point with positive coordinates.

\section{Conclusion}
\label{sec:future}

Despite the high computational cost associated to many of computations required to determine asymptotics using the methods of ACSV, the continued development of efficient computer algebra packages in Julia and other languages has made it feasible to automate the analysis beyond the simplest cases. There are many natural extensions still to be made, perhaps chiefly among them extending to the non-smooth case by incorporating algorithms for the Whitney stratification of algebraic varieties. Other interesting avenues for exploration include the development of better start systems for homotopy computations, to better match the number of critical points, and a theoretical study of the solution variety and its irreducible components for the monodromy approach (which could help the monodromy approach be competitive with or even surpass the polyhedral homotopy approach).

Finally, as already mentioned, in the combinatorial case the precision required to certify all minimal critical points with the same coordinate-wise modulus may not be practical. Our algorithm returns a warning with its output when it cannot verify equality between moduli to the precision required for rigorous certification.

\section*{Acknowledgments}
KL and SM acknowledge the support of the AMS Math Research Community \emph{Combinatorial Applications of Computational Geometry and Algebraic Topology}, which was funded by the National Science Foundation under Grant Number DMS 1641020. SM and JS's work partially supported by NSERC Discovery Grant RGPIN-2021-02382.
	
	\bibliographystyle{abbrv}
	\bibliography{references}
\end{document}